\documentclass[12pt,draftcls,onecolumn,letter]{IEEEtran}
\IEEEoverridecommandlockouts

\usepackage{lineno}
\usepackage{graphicx}
\usepackage{epsfig}
\usepackage{subfig}
\usepackage{times}
\usepackage{amsmath}
\usepackage{amssymb}
\usepackage{hyperref}

\usepackage{cite}
\usepackage{color}

\newtheorem{prob}{Problem}[section]

\title{Energy Generation and Distribution via Distributed Coordination: Case Studies}

\author{\small Hyo-Sung Ahn$^\dag$ and Byeong-Yeon Kim$^\dag$ 
\thanks{\small $^\dag$School of Mechatronics, Gwangju Institute of Science and Technology, Gwangju, South Korea.
{E-mail: hyosung@gist.ac.kr}}}

\begin{document}

\maketitle


\begin{abstract}
This paper presents case studies of the algorithms called ``energy generation and distribution via distributed coordination'', which was proposed in \cite{bykim_mesa_2012,bykim_dissertation_2013,bykim_arxiv_2014,bykim_submitted_2014}. For a convenience, we call ``energy generation and distribution via distributed coordination'' as ``DisCoord algorithms''. After concisely summarizing the ``DisCoord Algorithms'' in compact forms, we will list several scenarios for numerical tests. Then, while analyzing the simulation results, we will discuss the capability of the algorithms in handling the various cases. 
\end{abstract}

\section{Introduction}  \label{Introduction}
\begin{figure}[t]
\centering \includegraphics [scale=1.2]{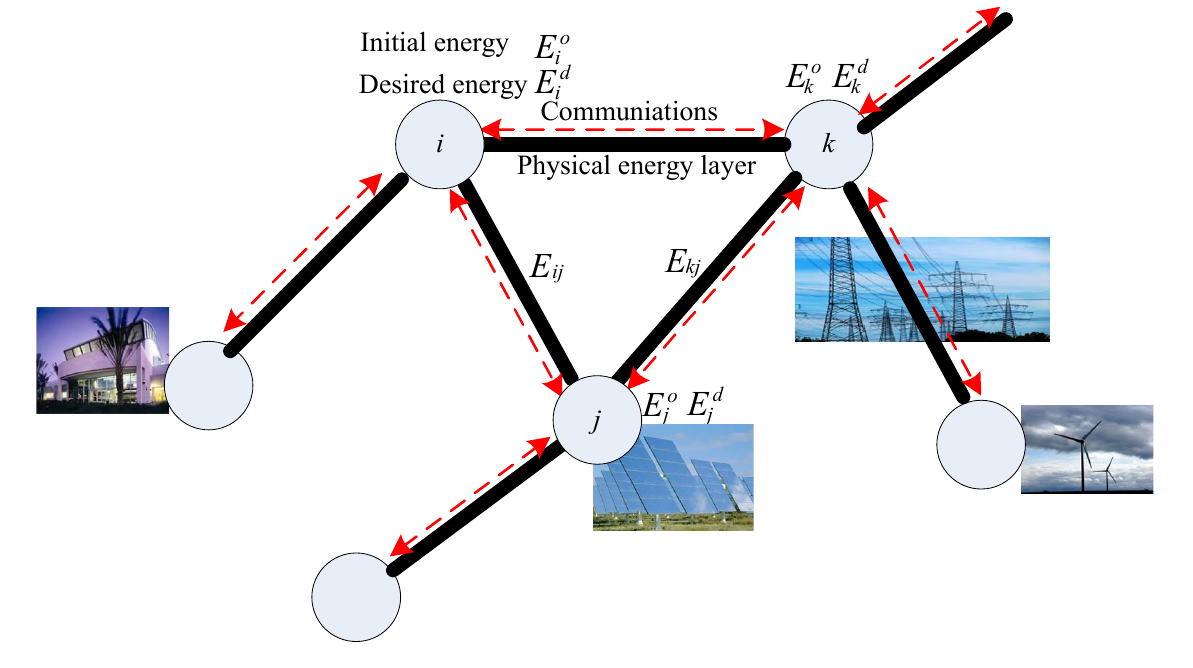}
\caption{Illustrative concept of energy generation and distribution in interconnected energy networks: Nodes can only communicate with their neighboring nodes to decide the amount of energy generation within some constrained ranges. After computing the amount of energy they have to generate, they also decide how to distribute energies to neighboring nodes via local interactions. The decision made on the basis of local interactions is called `distributed coordination'.} \label{energy_network}
\end{figure}

The DisCoord algorithms proposed in \cite{bykim_mesa_2012,bykim_dissertation_2013,bykim_arxiv_2014,bykim_submitted_2014} seek to solve a problem that is described in Fig.~\ref{energy_network}. The system of interest is composed of nodes distributed in space, communication network, and physical energy network. The nodes can generate a certain amount of energy within the node, while it has an initial energy in a certain level, $E_i^o$. The nodes also have some desired levels of energy, $E_i^d$. The nodes are connected through communication layer as well as through physical energy layer. Through the communication layer, they can exchange information about the status of individual node. Through the physical layer, the energy flows from a node to another node, $E_{ij}$. The energy flow is called energy distribution. The main purpose of DisCoord algorithms is to generate and distribute energies in an attempt to make each node to achieve the desired energy. Here, as the key constraints, there are the lower and upper levels in the generation of energy; also all the decision is made only through local interactions between neighboring nodes. This paper provides case studies taking account of various scenarios; while analyzing the simulation results, we will discuss the capability of the DisCoord algorithms. 

\section{DisCoord Algorithms}  \label{algorithm}
In DisCoord Algorithms, there are $N_{node}$ nodes, and each node is able to generate energy within some specified lower and upper boundaries. The nodes are interconnected for energy exchanges and for communications each other. To make the problem clear, the following symbols are provided.

\begin{tabular}{llll}
 $P_i$ & $i$-th node \\
 $E_i$ & the energy of $i$-th node \\ 
 $N_{node}$ & number of nodes \\ 
 $N_{edge}$ & number of edges \\   
 ${\mathcal N}$ & the set of nodes \\    
 ${\mathcal E}$ & the set of edges \\  
 ${\bf N}_i$ & the set of neighboring nodes of the $i$-th node \\
 $E_i^o$ & the initial energy of $i$-th node  \\   
 $E_i^d$ & the desired energy of $i$-th node  \\            
 $\Delta E_i$ & the energy generated by $i$-th node  \\ 
 $\underline{\Delta E_i}$ & the lower boundary over which $i$-th node has to generate\\  
 $\overline{\Delta E_i}$ & the upper boundary under which $i$-th node has to generate\\   
 $E_{ij}$ & energy flow from the $i$-th node to the $j$-th node  \\                 
\end{tabular}

\ \

The combination of the set of nodes and set of edges is called a graph, i.e., ${\mathcal G}\doteq \{{\mathcal N}, {\mathcal E}\}$. Thus, the overall topology of the system is defined by the graph ${\mathcal G}$. 

The DisCoord Algorithms seek to solve the following problem \cite{bykim_mesa_2012}.

\begin{prob}\label{prob1}
Given initial and desired energies of each nodes, $E_i^o$ and $E_i^d$, generate energies $\Delta E_i$ and exchange energies, $E_{ji}$, with neighboring nodes such that $E_i^o \rightarrow E_i^d$, by relying upon only local interactions. 
\end{prob}

In the above problem, the energy generation should be done with the constraints of 
\begin{align} \label{generation_constraint}
\underline{\Delta E_i} \leq \Delta E_i  \leq \overline{\Delta E_i}
\end{align} 

For the energy generation, the node $i$ can exchange a certain type of information only with its neighboring nodes. It is supposed that the information is exchanged through communication layer network (or cyber-layer network) \cite{bykim_submitted_2014}. The energy exchange means the energy flow between the neighboring nodes; so the energy exchange is done in physical-layer network. 

It is also important to notice that there is a relationship between the desired energies and initial and generated energies. That is, the sum of the desired energies should be less than the sum of initial energies and the maximum energy that can be generated by the nodes. It can be constrained as follows:
\begin{align} \label{desired_constraint}
\sum_{i=1}^{N_{node}} ( E_i^o + \underline{\Delta E_i})  \leq \sum_{i=1}^{N_{node}} E_i^d \leq \sum_{i=1}^{N_{node}} ( E_i^o + \overline{\Delta E_i})
\end{align}

Eventually, the DisCoord Algorithms should generate the energies ${\Delta E_i}$ such as 
\begin{align} \label{supply_demand_balance}
\sum_{i=1}^{N_{node}} (E_i^o + {\Delta E_i}) = \sum_{i=1}^{N_{node}} E_i^d 
\end{align}
which is called the supply-demand balance. So, given the desired energies of individual nodes, one of the key constraints in the DisCoord Algorithms is to satisfy the supply-demand balance during the generations of ${\Delta E_i}$.

\subsection{Energy generation}
The purpose of energy generation is to generate ${\Delta E_i}$ only using local neighboring interactions under the constraints of (\ref{generation_constraint}) and (\ref{supply_demand_balance}). The algorithm is composed of three steps. In the following algorithm, the interim parameters $z_i$ and $w_i$ are used. 

\begin{enumerate}
\item Initial values \\
-- $z_i(0)= E_i^d  -E_i^o -\underline{\Delta E_i} $ \\
-- $w_i(0) = \overline{\Delta E_i} -\underline{\Delta E_i}$ 

\item Find the steady state solutions from the following dynamic equations \\
-- $z_i(t + 1) = \frac{1}{1 + |{{\bf N}_i}|}  z_i(t) + \sum\limits_{j \in {{\bf N}_i}}
\frac{1}{1 +|{{\bf N}_j}|} z_j(t) $        \\
-- $w_i(t + 1) = \frac{1}{1 + |{{\bf N}_i}|}  w_i(t) + \sum\limits_{j \in {{\bf N}_i}}
\frac{1}{1 +|{{\bf N}_j}|} w_j(t) $ \\
where $|{{\bf N}_i}|$ is the cardinality of the set ${{\bf N}_i}$. Let the steady-state solutions of the above equations denote as $z_i(\infty) \doteq \lim_{t \rightarrow \infty} z_i(t)$ and $w_i(\infty) \doteq \lim_{t \rightarrow \infty} w_i(t)$.

\item Energy generation of nodes \\
-- $\Delta E_i = \underline{\Delta E_i} + \frac{(\overline{\Delta E_i} - \underline{\Delta E_i})}{w_i(\infty)} z_i(\infty)$
\end{enumerate}

Note that in the above $2$nd step, for the update of $z_i(t + 1)$ and $w_i(t + 1)$, the node $i$ needs information of $z_j(t)$ and $w_j(t)$, which are the values of neighboring interim parameters. So, to compute $\Delta E_i$ at the node $i$, it needs to communicate with the neighboring nodes. Since $|{{\bf N}_i}|$ and $|{{\bf N}_j}|$ are fixed and can be exchanged by initial local interactions, it is simply assumed that they are available to neighboring nodes.

\subsection{Energy distribution}
The purpose of the energy distribution is to make the level of energy of each node become $E_i^d$, i.e., $E_i \rightarrow E_i^d$ by the energy flows $E_{ji}$. This algorithm is also composed of three steps. In the following algorithm, the interim parameters $h_{ij}$ and $g_{i}$ are used. 

\begin{enumerate}
\item Initial values \\
-- $h_{ij}(0)= 0$ \\
-- $g_i(0) = E_i^o  +   \Delta E_i - E_i^d$ \\
where $\Delta E_i$ is computed from the $3$-rd step in the energy generation algorithm. 
\item Find the steady state solutions from the following dynamic equations \\
-- $h_{ij}(t + 1) = h_{ij}(t) + a_{ij}(g_j(t) -g_i(t))$ \\
-- $g_i(t + 1) = g_i(t) + \sum\limits_{j \in {\bf N}_i} a_{ij}(g_j(t) - g_i(t))$ \\
where $a_{ij} = \frac{1} {1 + \max \left\{ | {{\bf N}_i} |, |{\bf N}_j | \right\}}$. Let the steady-state solution of the above equation denote $h_{ij}(\infty) \doteq \lim_{t \rightarrow \infty} h_{ij}(t)$. 
\item Energy flow from $i$-th node to the $j$-th node \\
-- $E_{ij} =  h_{ij}(\infty)$ 
\end{enumerate}

Note that in the above $2$nd step, for the update of $h_{ij}(t + 1)$ and $g_i(t + 1)$, the node $i$ needs information of $g_j(t)$ from the neighboring nodes. Since $a_{ij}$ are fixed and can be also exchanged by local interactions initially, it is assumed that they are available to neighboring nodes.

\subsection{Comments on algorithms}
The convergence of algorithms has been completely proved in \cite{bykim_dissertation_2013,bykim_arxiv_2014,bykim_submitted_2014}. It is assumed that the desired energy of individual node is given as (\ref{desired_constraint}). The energy generation algorithm computes $\Delta E_i$; and then $\Delta E_i$ is used for the computation of $E_{ij}$ at the energy distribution algorithm. It is also noticeable that the solutions of the energy generation and energy distribution algorithms are computed through communications among neighboring nodes. So, after obtaining the solutions $\Delta E_i$ and $E_{ij}$, each node generates physical energy and distributes energy to the neighboring nodes simultaneously.

\section{Case studies}

\subsection{Case -1: Balanced case} In this case, the sum of the desired energies of all nodes is well given as 
$\sum_{i=1}^{N_{node}} ( E_i^o + \underline{\Delta E_i})  \leq \sum_{i=1}^{N_{node}} E_i^d \leq \sum_{i=1}^{N_{node}} ( E_i^o + \overline{\Delta E_i})$.  For the balanced case, we consider the following scenario:\\
-- Six nodes: $P_1, P_2, P_3, P_4, P_5$, and $P_6$ \\
-- Seven undirected edges: $\mathcal{E} \doteq \{ e(1,2), e(2,3), e(3,4), e(1,5), e(3,5), e(4,5), e(4,6)\}$ \\
-- Initial energies: $E_1^o=0, E_2^o=10, E_3^o=10, E_4^o=2, E_5^o=0$, and $E_6^o=10$ \\ 
-- Desired energies: $E_1^d=5, E_2^d=15, E_3^d=20, E_4^d=30, E_5^d=2$, and $E_6^d=20$ \\ 
-- Lower boundaries of generation capability:  $\underline{\Delta E_i}=0$ for all $i$ \\
-- Upper boundaries of generation capability:  $\overline{\Delta E_1}=5, \overline{\Delta E_2}=15, \overline{\Delta E_3}=15, \overline{\Delta E_4}=15, \overline{\Delta E_5}=20$ and $\overline{\Delta E_6}=15$ 

\begin{figure}[t]
\centering \includegraphics [scale=0.7]{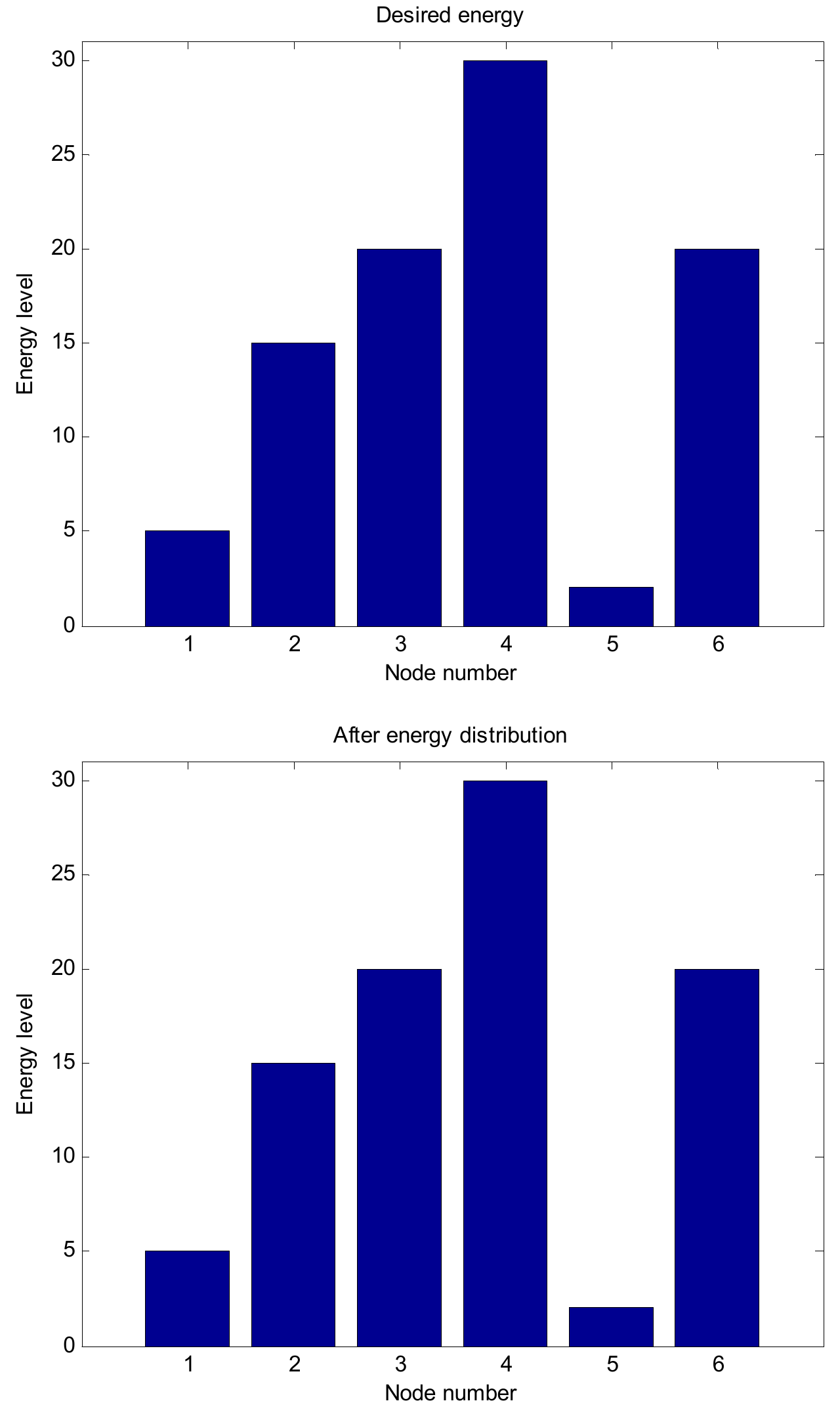}
\caption{Case $1$: Desired energies(upper figure) and achieved energies(lower figure).} \label{case11}
\end{figure}

\begin{figure}[t]
\centering \includegraphics [scale=0.7]{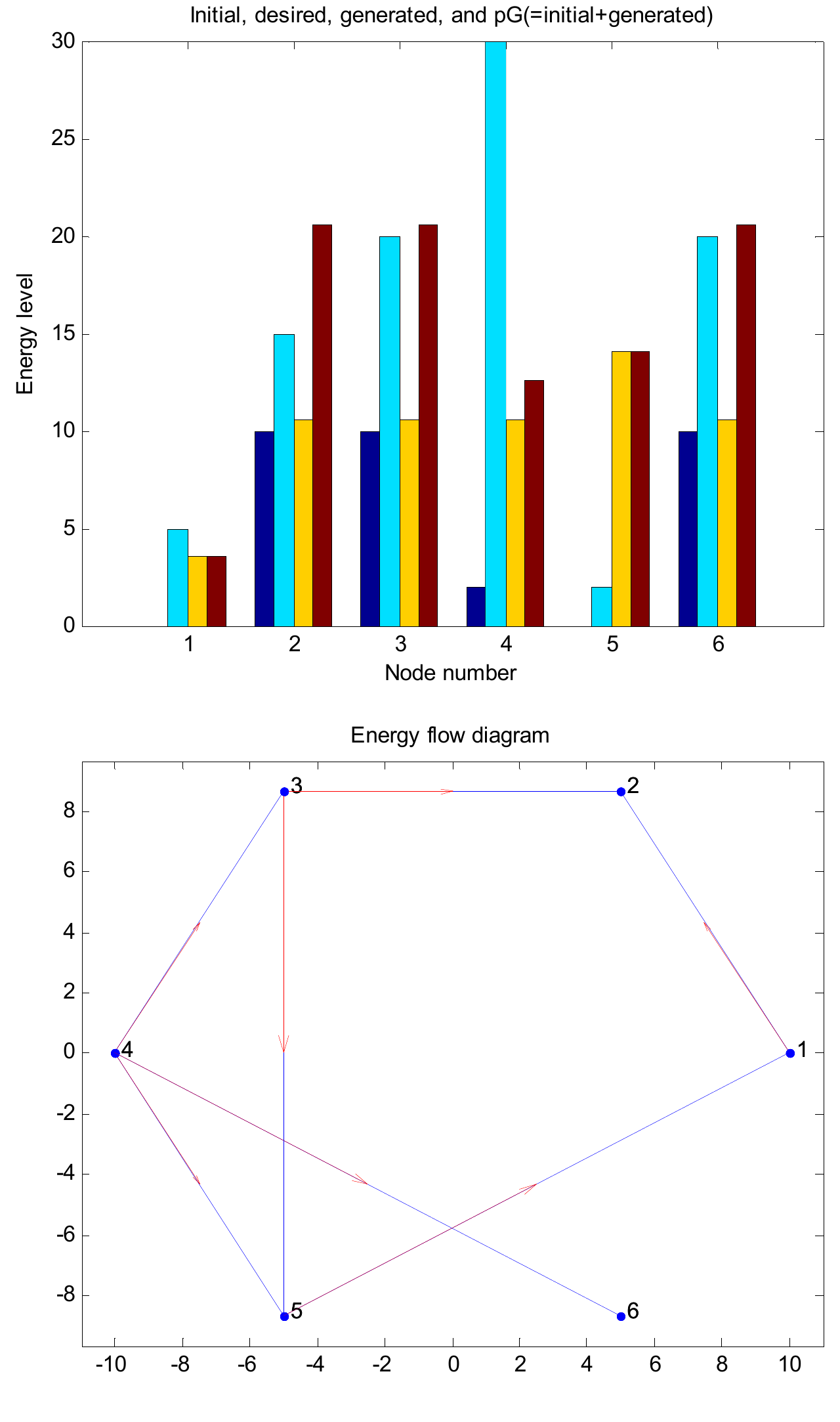}
\caption{Case $1$: Upper figure - Initial(blue), desired(aqua), generated(yellow), and initial$+$generated(red); Lower figure - Energy flow diagram.} \label{case12}
\end{figure}

Clearly, since $\sum_{i=1}^{N_{node}} ( E_i^o + \underline{\Delta E_i}) =32$, $\sum_{i=1}^{N_{node}} E_i^d = 92$, and $\sum_{i=1}^{N_{node}} ( E_i^o + \overline{\Delta E_i}) =117$, the supply-demand balance is satisfied. Now using the DisCoord algorithms given in the previous section, we obtain\\
-- The generated energies: $\Delta E_1=3.5294$, $\Delta E_2=10.5882$, $\Delta E_3=10.5882$, $\Delta E_4=  10.5882$, $\Delta E_5=14.1176$, $\Delta E_6=10.5882$  \\ 
-- Energy flows:   $E_{12}=1.6298$,  $E_{32}=3.9584$,  $E_{35}= 2.5768$,  $E_{43}=7.1234$,   $E_{45}=9.7001$, $E_{46}=0.5882$, $E_{51}=0.1593$     \\                 
From the above result, we see that the total amount of energy flows is $25.7360$. As shown in Fig.~\ref{case11}, the desired energies have been well achieved after energy distribution. Fig.~\ref{case12} shows the overall variation in energy level of each node, and overall energy flows in the graph. 

\subsection{Case -2: Under-demand case} 
\begin{figure}[t]
\centering \includegraphics [scale=0.7]{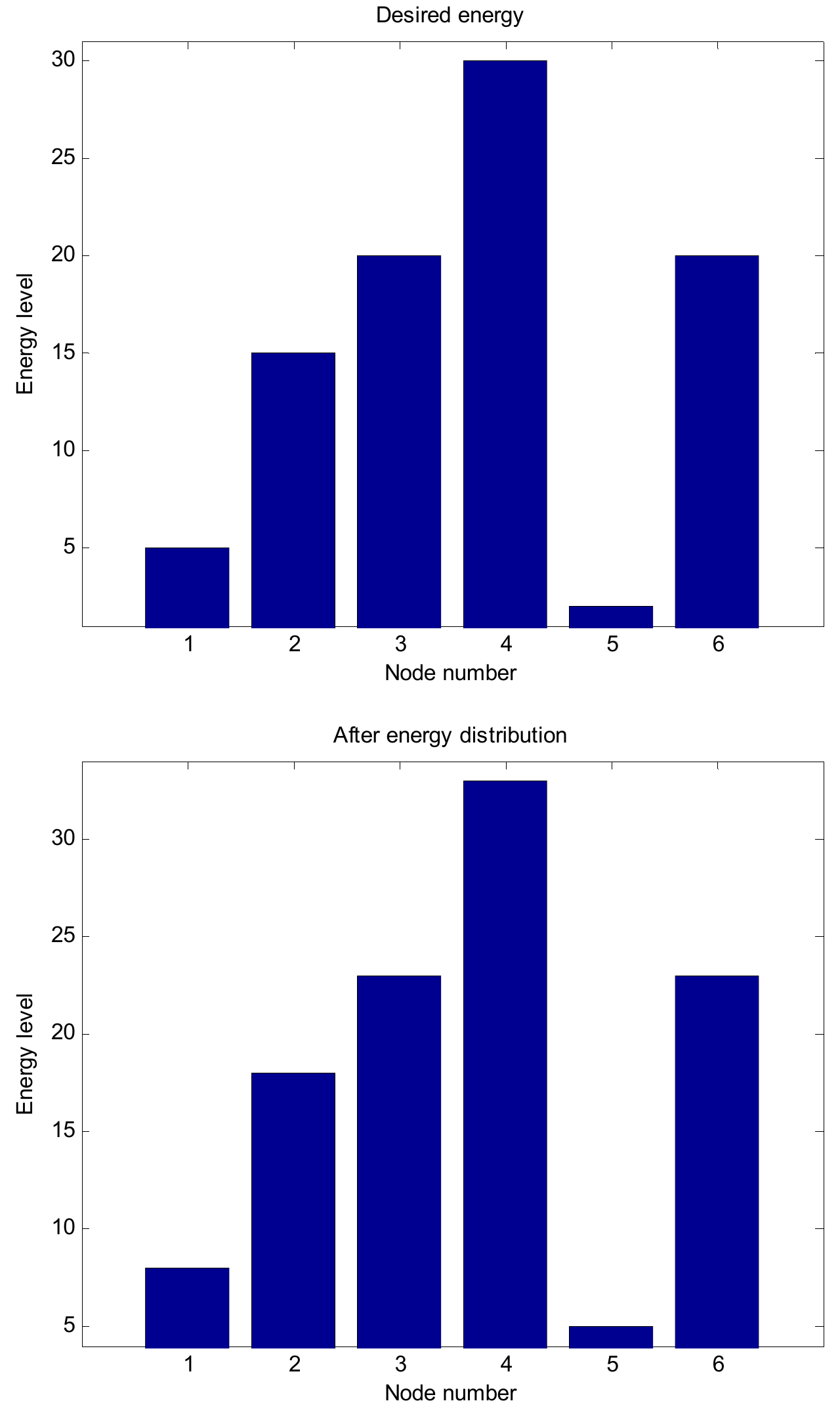}
\caption{Case $2$: Desired energies(upper figure) and achieved energies(lower figure).} \label{case21}
\end{figure}

\begin{figure}[t]
\centering \includegraphics [scale=0.7]{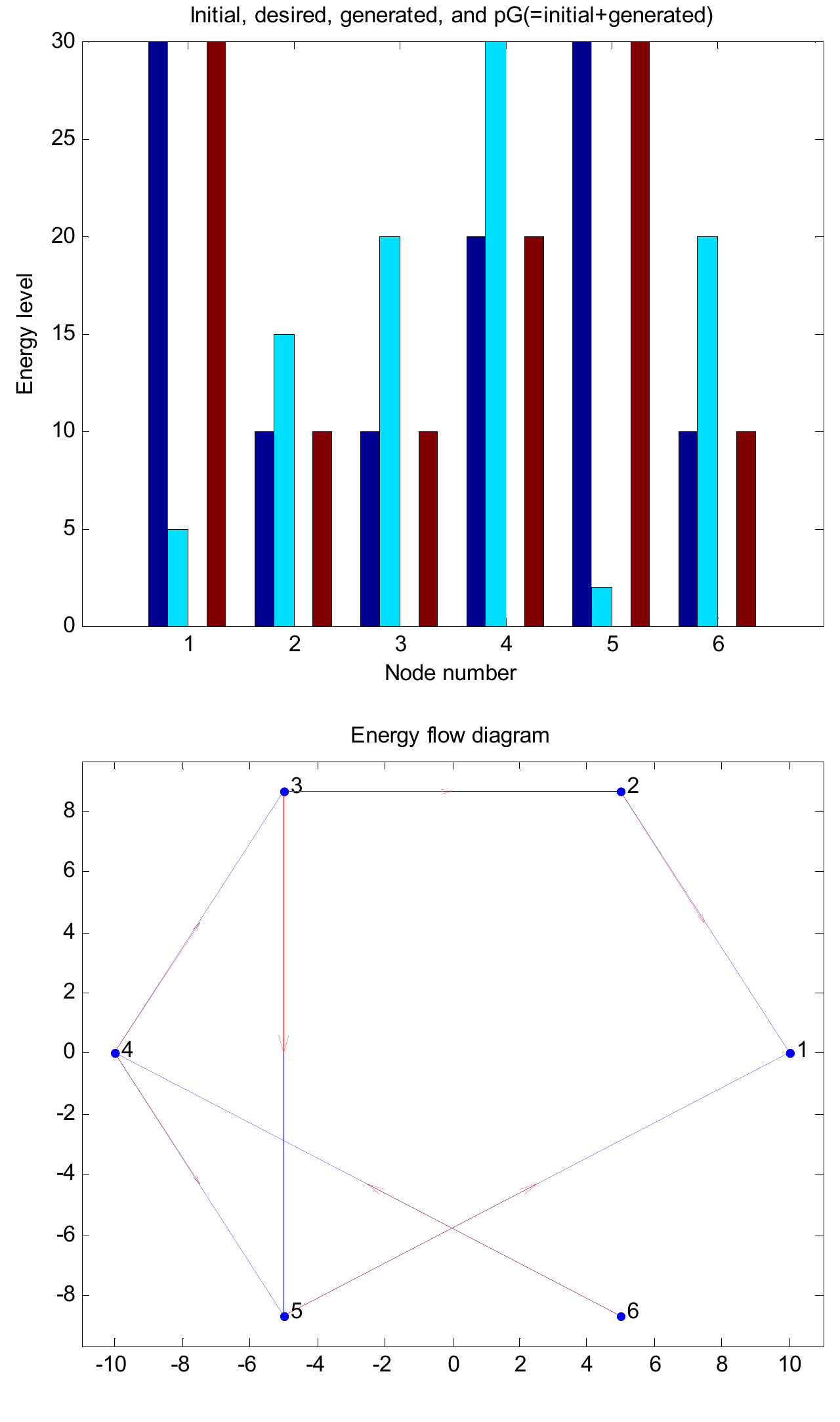}
\caption{Case $2$: Upper figure - Initial(blue), desired(aqua), generated(yellow), and initial$+$generated(red); Lower figure - Energy flow diagram.} \label{case22}
\end{figure}

In this case, the sum of the desired energy is less than the minimum energy boundary, i.e., $\sum_{i=1}^{N_{node}} E_i^d < \sum_{i=1}^{N_{node}} ( E_i^o + \underline{\Delta E_i})$. In the simulation, the same scenario as the case -$1$ is considered except the initial energies \\
-- Initial energies: $E_1^o=30, E_2^o=10, E_3^o=10, E_4^o=20, E_5^o=30$, and $E_6^o=10$ \\ 
So, we have $\sum_{i=1}^{N_{node}} E_i^d = 92 < \sum_{i=1}^{N_{node}} ( E_i^o + \underline{\Delta E_i}) = 110$. From the DisCoord algorithms, we obtain\\
-- The generated energies: $\Delta E_1=0$, $\Delta E_2=0$, $\Delta E_3=0$, $\Delta E_4=0$, $\Delta E_5=0$, $\Delta E_6=0$  \\ 
-- Energy flows:   $E_{21}=15.4146$,  $E_{32}=7.4146$,  $E_{35}= 12.3902$,  $E_{43}=6.8049$,   $E_{45}=19.1951$, $E_{51}=6.5854$, $E_{64}=13.0000$     \\                 
From the above result, we see that the total amount of energy flows is $80.8049$. After the energy distribution, it is observed that the errors, i.e., $e_i = E_i^d - E_i$, are $e_1 = -3, e_2=-3, e_3 = -3, e_4 = -3, e_5 = -3, e_6= -3$. Fig.~\ref{case21} shows that the desired energies cannot be achieved because there are over energies in the network. So, as shown in the lower figure of Fig.~\ref{case21}, after the energy distribution, each node still  has more energy than the desired energy level. It is interesting to observe that the errors of each node are same as $-3$; so it seems that the DisCoord Algorithms evenly distribute the remaining energies to the network. As shown in Fig.~\ref{case22}, clearly there is no energy generation from the DisCoord Algorithms. Thus, we can see that when there are over energy in the network, the algorithms do not generate energy in the whole network. 

\subsection{Case -3: Over-demand case} 
\begin{figure}[t]
\centering \includegraphics [scale=0.7]{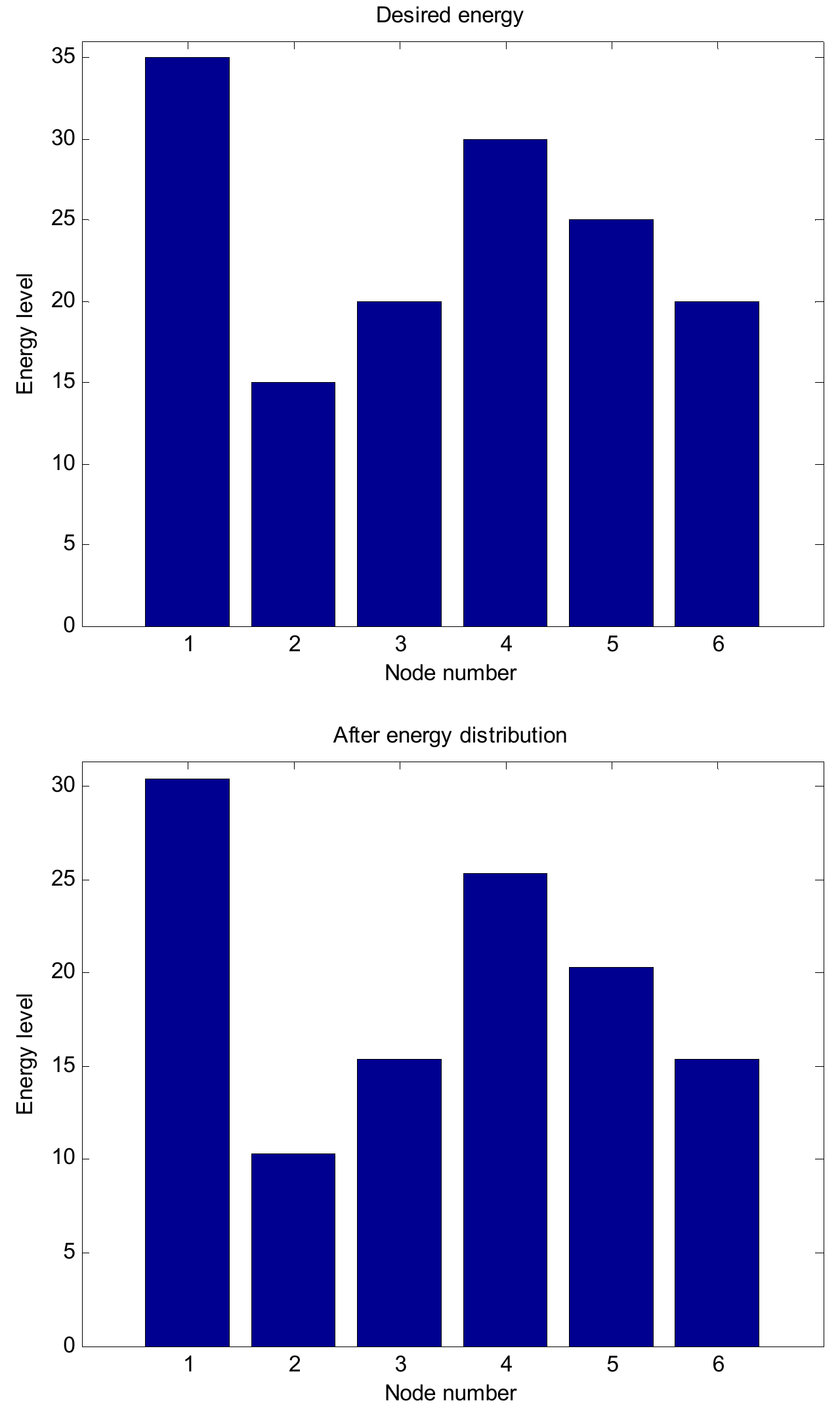}
\caption{Case $3$: Desired energies(upper figure) and achieved energies(lower figure).} \label{case31}
\end{figure}

\begin{figure}[t]
\centering \includegraphics [scale=0.7]{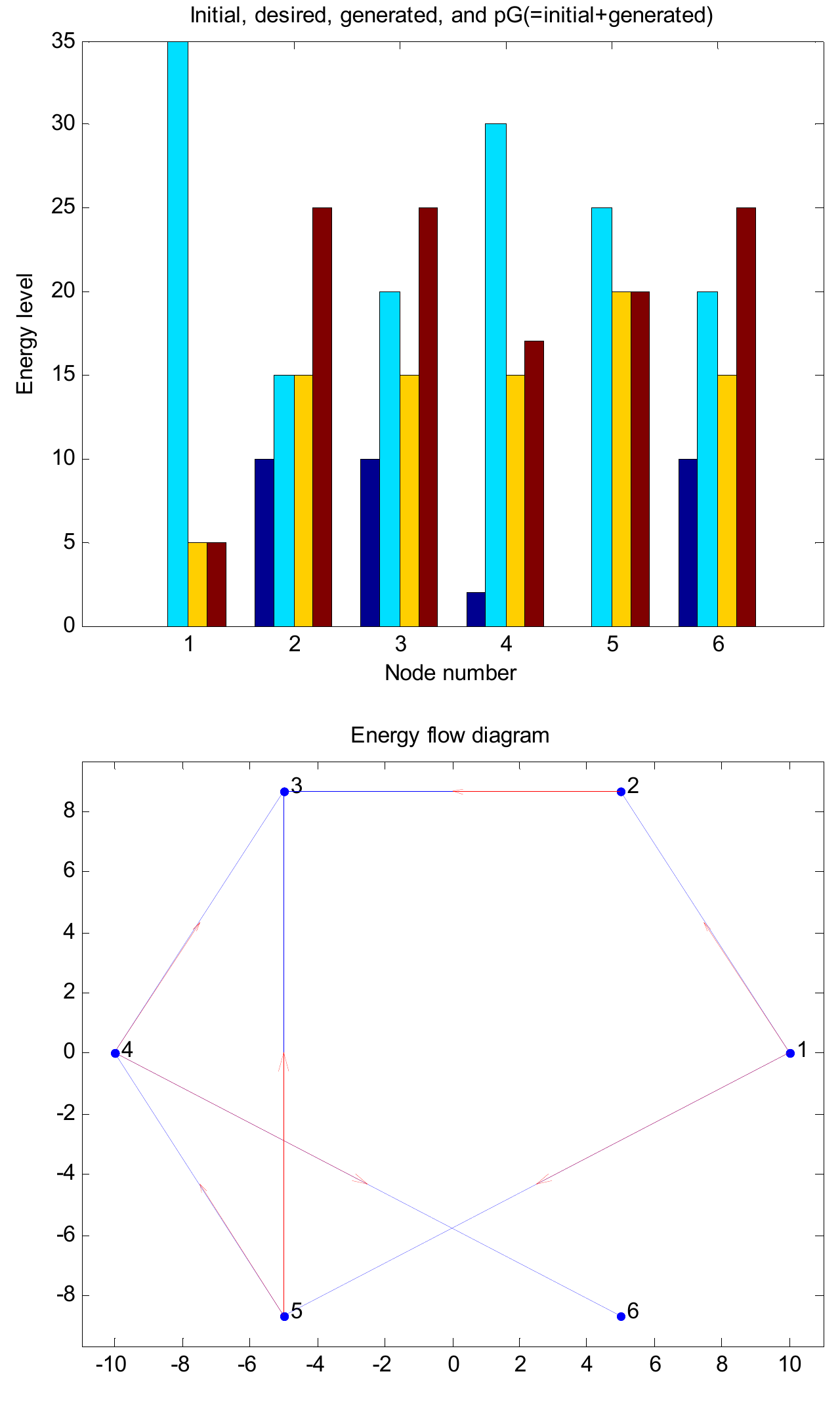}
\caption{Case $3$: Upper figure - Initial(blue), desired(aqua), generated(yellow), and initial$+$generated(red); Lower figure - Energy flow diagram.} \label{case32}
\end{figure}

In this case, the sum of the desired energy is greater than the maximum energy boundary, i.e., 
$\sum_{i=1}^{N_{node}} ( E_i^o + \overline{\Delta E_i})< \sum_{i=1}^{N_{node}} E_i^d$.  
For the simulation, the same scenario as the case -$1$ is considered except the desired energies \\
-- Desired energies: $E_1^d=35, E_2^d=15, E_3^d=20, E_4^d=30, E_5^d=25$, and $E_6^d=20$ \\ 
So, we have $\sum_{i=1}^{N_{node}} ( E_i^o + \overline{\Delta E_i}) = 117 < \sum_{i=1}^{N_{node}} E_i^d = 145$. From the DisCoord algorithms, we obtain\\
-- The generated energies: $\Delta E_1=5$, $\Delta E_2=15$, $\Delta E_3=15$, $\Delta E_4=15$, $\Delta E_5=20$, $\Delta E_6=15$  \\ 
-- Energy flow:   $E_{12}=16.5854$,  $E_{15}=8.7480$,  $E_{23}= 1.9187$,  $E_{43}= 2.1382$,   $E_{46}=9.6667$, $E_{53}=5.6098$, $E_{54}=3.4715$  \\
From the above result, we see that the total amount of energy flows is $48.1382$. After the energy distribution, it is observed that the errors are $e_1 = 4.6667, e_2= 4.6667, e_3 = 4.6667, e_4 = 4.6667, e_5 =  4.6667, e_6= 4.6667$. Fig.~\ref{case31} shows that the desired energies at each node are not achieved due to the energy deficiency. Fig.~\ref{case32} reveals that the energy at each node is generated in maximum. Similarly to the Case $2$, the DisCoord Algorithms still attempt to distribute the energies to the network evenly. So, the errors at each node are same. The algorithms have attempted to generate the maximum energy; after generating the maximum energy, the energies have been distributed evenly to the network to make the same errors at all nodes.

\subsection{Case -4: No generation case} In this case, it is supposed that 
$\sum_{i=1}^{N_{node}}  E_i^o = \sum_{i=1}^{N_{node}} E_i^d$. Then, it does not need to generate energy; only by distributing the initial energy among the nodes, the desired energy of individual node can be achieved. For the simulation, the same scenario as the case -$1$ is considered except the desired energies \\
-- Initial energies: $E_1^o=30, E_2^o=10, E_3^o=10, E_4^o=20, E_5^o=12$, and $E_6^o=10$ \\ 
So, the sum of initial energies is equal to the sum of the desired energies. From the DisCoord algorithms, we obtain\\
-- The generated energies are all zero. \\
-- Energy flows:   $E_{21}=13.6585$,  $E_{32}=8.6585$,  $E_{35}= 7.5610$,  $E_{43}= 6.2195$,   $E_{45}= 13.7805$, $E_{51}= 11.3415$, $E_{64}=10$  \\
From the above result, we see that the total amount of energy flows is $71.2195$. 
\begin{figure}[t]
\centering \includegraphics [scale=0.7]{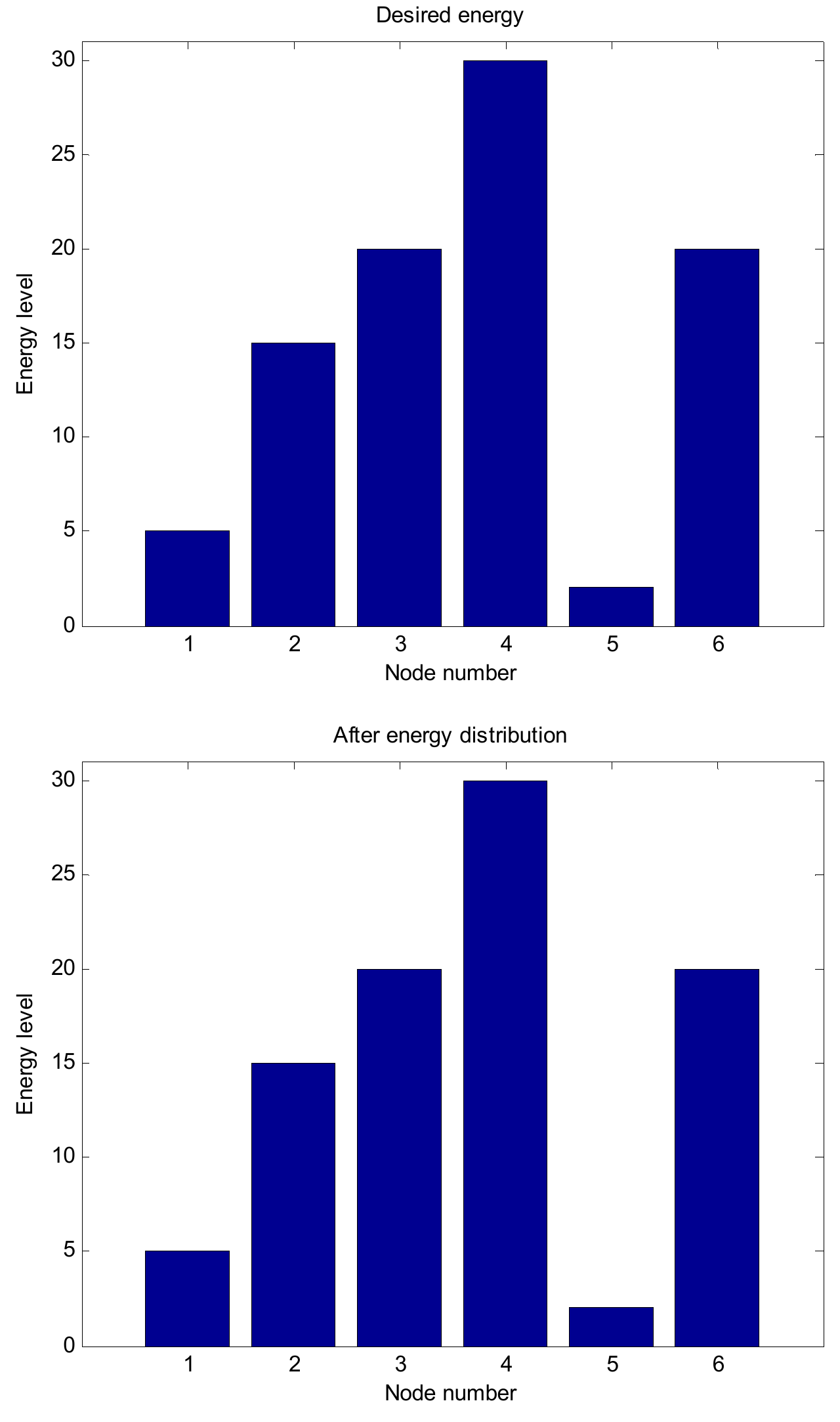}
\caption{Case $4$: Desired energies(upper figure) and achieved energies(lower figure).} \label{case41}
\end{figure}
\begin{figure}[t]
\centering \includegraphics [scale=0.7]{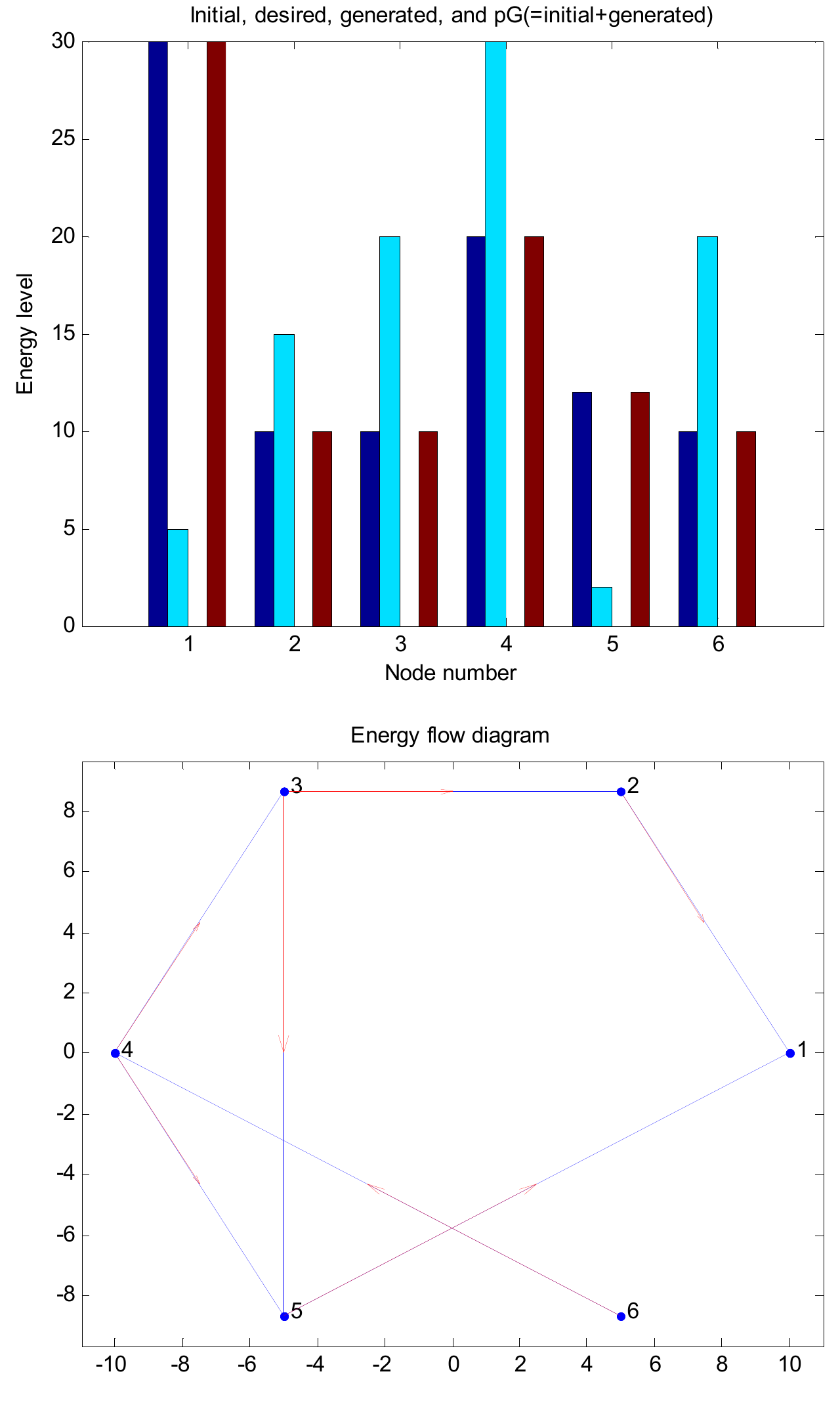}
\caption{Case $4$: Upper figure - Initial(blue), desired(aqua), generated(yellow), and initial$+$generated(red); Lower figure - Energy flow diagram.} \label{case42}
\end{figure}
Fig.~\ref{case41} shows that the desired energies have been well achieved without generating any energy. So, as expected, the DisCoord Algorithms only distribute the initial energies to the network to make errors of each node zero.  Fig.~\ref{case42} also shows that there is no energy generation; but only there are energy flows. 

\section{Conclusion}  \label{conclusion}
This paper has presented several cases as applications of DisCoord algorithms. As shown in the previous examples, the DisCoord algorithms can deal with various cases even though the supply-demand balance is not ensured any more. The main capability of the algorithms is to decide all the generation and distribution through local interactions, i.e., by distributed coordination.

\section{Acknowledgement}
It is recommended to see `{Byeong-Yeon Kim, ``Coordination and control for energy distribution using consensus algorithms in interconnected grid networks'', Ph.D. Dissertation, School of Information and Mechatronics, Gwangju Institute of Science and Technology, 2013} \cite{bykim_dissertation_2013}' for applications to various engineering problems of the algorithms developed in this paper. The reader can send email to hyosung@gist.ac.kr to get the MATLAB source code.



\begin{thebibliography}{99}
\bibitem{bykim_mesa_2012}
B.-Y. Kim, K.-K. Oh, and H.-S. Ahn,
\newblock ``Power Distribution with Consensus,''
\newblock {\em Proc. of the 2012 IEEE/ASME Int. Conf. Mechatronics and Embedded Systems and Applications}, Suzhou, China, July 8-10, 2012.

\bibitem{bykim_dissertation_2013}
B.-Y. Kim,
\newblock ``Coordination and control for energy distribution using consensus algorithms in interconnected grid networks,'' 
\newblock in {\em Ph.D Dissertation}, School of Information and Mechatronics, Gwangju Institute of Science and Technology, 2013.

\bibitem{bykim_arxiv_2014}
B.-Y. Kim, K.-K. Oh, and H.-S. Ahn,
\newblock ``Power Generation and Distribution via Distributed Coordination Control,''
\newblock {\em arXiv:1407.4870 [math.OC]}, 2014.

\bibitem{bykim_submitted_2014}
B.-Y. Kim, K.-K. Oh, and H.-S. Ahn,
\newblock ``Coordination and Control for Energy Distribution in Distributed Grid Networks,''
\newblock {\em submitted for a publication}, 2014.
\end{thebibliography}

\end{document}